\documentstyle[12pt]{article}
\newtheorem{Th}{Theorem}
\newtheorem{Prop}[Th]{Proposition}

\newtheorem{Cor}[Th]{Corollary}
\newtheorem{Lemma}[Th]{Lemma}
\newcommand{\Proof}{\noindent {\it Proof. }}

\newcommand{\be}{\begin{eqnarray*}}
\newcommand{\ee}{\end{eqnarray*}}
\newcommand{\pkEF}{{\cal P}(^k\!E,F)}

\newcommand{\pcckEF}{{\cal P}_{cc}(^k\!E,F)}
\newcommand{\pcck}{{\cal P}_{cc}(^k\!}
\newcommand{\puckEF}{{\cal P}_{uc}(^k\!E,F)}
\newcommand{\LkEF}{{\cal L}(^k\!E,F)}

\newcommand{\WCo}{{\cal WC}o\,}
\newcommand{\CC}{{\cal CC}}
\newcommand{\UC}{{\cal UC}}
\newcommand{\LEF}{{\cal L}(E,F)}

\newcommand{\finesp}{\hspace*{\fill} $\Box$\vspace{0.5\baselineskip}}
\newcommand{\fin}{\hspace*{\fill} $\Box$}

\newcommand{\N}{\bf N}
\newcommand{\lra}{\longrightarrow}

\newcommand{\ra}{\rightarrow}

\newcommand{\bmi}{\noindent \begin{minipage}[t]{136mm}}

\newcommand{\emi}{\end{minipage}}

\begin{document}

\title{When every polynomial is \\
      unconditionally converging}

\author{Manuel Gonz\'alez\thanks{Supported in part by DGICYT Grant PB
      91--0307 (Spain)} \\
      Departamento de Matem\'aticas \\
      Facultad de Ciencias\\
      Universidad de Cantabria \\ 39071 Santander (Spain)
       \and
      Joaqu\'\i n  M. Guti\'errez\thanks{Supported in part by DGICYT
      Grants PB 90--0044 and PB 91-0307 (Spain)}\\
      Departamento de Matem\'atica Aplicada\\
      ETS de Ingenieros Industriales \\
      Universidad Polit\'ecnica de Madrid\\
      C. Jos\'e Guti\'errez Abascal 2 \\
      28006 Madrid (Spain)}

\date{\mbox{}}
\maketitle

\begin{abstract}
Letting $E$, $F$ be Banach spaces,
the main two results of this paper are the following:
(1) If every (linear bounded) operator $E\ra F$ is unconditionally converging,
then every polynomial from $E$ to $F$ is unconditionally converging (definition
as in the linear case).
(2) If $E$ has the Dunford-Pettis property and every operator $E\ra F$ is weakly
compact, then every $k$-linear mapping from $E^k$ into $F$ takes weak Cauchy
sequences into norm convergent sequences. In particular, every polynomial from
$\ell_\infty$ into a space containing no copy of $\ell_\infty$ is completely
continuous.
This solves a problem raised by the authors in a previous paper, where they
showed that there exist nonweakly compact polynomials from $\ell_\infty$ into
any nonreflexive space.
\end{abstract}

\section{Introduction}
\label{Intr}

Throughout, $E$, $F$ will be Banach spaces.
We denote by $\LEF$ the space of all (linear bounded) operators from $E$
to $F$, and by $\WCo (E,F)$ the subspace of all weakly compact operators.
We say that $T \in \LEF$ is {\it completely continuous\/} if it takes weakly
convergent sequences into norm convergent sequences, and $T$ is {\it
unconditionally converging\/} if it takes weakly unconditionally Cauchy (w.u.C.)
series into unconditionally convergent (u.c.) series
(the definitions are recalled below).
The respective subspaces of operators are denoted $\CC (E,F)$ and $\UC (E,F)$.

It is well known that if $F$ contains no copy of $\ell_\infty$, then we have the
equalities
$$
{\cal L}(\ell_\infty,F) = \WCo (\ell_\infty,F) = \CC (\ell_\infty,F)
= \UC (\ell_\infty,F) \, .
$$

We denote by $\pkEF$ the space of all $k$-homogeneous (continuous) polynomials
from $E$ to $F$, and by $\pcckEF$ the subspace of completely continuous
polynomials, i.e., the polynomials taking weakly convergent sequences into norm
convergent ones.

In \cite{GG} the authors showed that, contrarily to the linear case $(k=1)$,
whenever
$F$ is nonreflexive, for every integer $k \geq 2$, there is a polynomial
$P\in {\cal P}_{cc}(^k\!\ell_\infty , F)$ which is not weakly compact.
It can be obtained as the composition of the following three mappings
$$
\ell_\infty \stackrel{U}{\lra} \ell_2 \stackrel{Q}{\lra} \ell_1
\stackrel{T}{\lra} F
$$
where $U$ is a completely continuous linear surjection, $Q$ is the polynomial
given by $Q\left( (x_n)_n\right) = \left( x^k_n\right) _n$, and $T$ is a
quotient onto a separable nonreflexive subspace of $F$.

It is then natural to ask the following question:

(a) is every polynomial from $\ell_\infty$ to $c_0$ completely continuous?

The authors also proved \cite{GG} that, given $P\in\pkEF$, if $\sum x_i$
is a w.u.C. (resp., u.c.) series in $E$, then $\sum P(x_i)$ is a w.u.C.
(resp., u.c.)
series in $F$.

This justifies the introduction of the following class of
polynomials: we say that $P$ is {\it unconditionally converging,} and write
$P\in\puckEF$, if for every w.u.C. series
$\sum x_i$ in $E$, the series $\sum P(x_i)$ is u.c.\ in $F$.
Several properties of Banach spaces are obtained in~\cite{GG} in terms of the
relationship of $\puckEF$ with other classes of polynomials used in the
literature.
It is proved in particular that we always have $\pcckEF \subseteq \puckEF$,
and that every weakly compact polynomial is unconditionally converging.

Therefore, previous to question (a) is the following:

(b) is every polynomial from $\ell_\infty$ to $c_0$
unconditionally converging?

In the present paper, we first prove (Section~\ref{PUC})
that whenever $E$ and $F$ satisfy the condition $\LEF = \UC (E,F)$,
we also have $\pkEF = \puckEF$ for all $k\in\N$, where $\N$ denotes the natural
numbers.

Recall that $E$ has the {\it Dunford-Pettis
property\/} (DPP) if for every $F$ we have
$\WCo (E,F) \subseteq \CC (E,F)$.
Our second main result (Section~\ref{PCC})
states that, whenever $E$ has the
DPP, and $\LEF = \WCo (E,F)$, we then have $\pkEF = \pcckEF$.
This happens, for instance, if $E=C(K)$ with $K$ stonean, and $F$ contains no
copy of $\ell_\infty$ $(F \not\supset\ell_\infty )$; also, if $E=C(K)$ and
$F\not\supset c_0$.
It may be useful to recall that the problem of finding pairs $E$, $F$ so that
$\pkEF =\pcckEF$ has received some attention. Pe\l czy\'nski showed~\cite{Pel1}
that this is the case for $E=\ell_p$ and $F=\ell_q$, with $kq<p$.
Gonzalo and Jaramillo \cite{GJ} have recently extended this result to spaces
admitting upper and lower $p$-estimates.

In Section \ref{exten},
we prove the existence of a completely continuous
extension of every polynomial $P\in \pcckEF$ to the
bidual of $E$, when the dual of $E$ has the DPP.
The proofs of all these results take advantage of an idea of~\cite{Ry}.

We denote by $E^*$ the dual of $E$, and by $T^*: F^* \ra E^*$ the adjoint of the
operator $T:E\ra F$. The space of $k$-linear (continuous)
mappings from $E^k$ into $F$ is denoted by
$\LkEF$.
To each $P\in \pkEF$ we can associate a unique symmetric $\hat{P}\in\LkEF$
so that $P(x) = \hat{P}(x,\ldots ,x)$ for all $x\in E$.
We say that $A\in \LkEF$ is {\it completely continuous\/} if given weak Cauchy
sequences $\left( x_1^n\right) , \ldots ,\left( x_k^n\right) \subset E$,
the sequence $\left( A(x_1^n,\ldots ,x_k^n)\right) _n$ is norm convergent in
$F$.

It is well known \cite[Proposition~5]{Pel2} that if $E$ has the DPP,
then every
$k$-linear mapping from $E^k$ to $F$ takes weak Cauchy sequences into
weak Cauchy sequences. This is not true in general: for instance, the
bilinear mapping
$$
\left( \left(x_n\right) , \left(y_n\right)\right) \in \ell_2
\mapsto \left(x_ny_n\right)\in \ell_1
$$
takes a weakly null sequence into a sequence having no weak Cauchy subsequence.

For the general theory of polynomials on Banach spaces,
we refer to \cite{Mu}.
Finally, let us recall that a formal series $\sum x_i$ in $E$ is w.u.C. if
for every $\phi\in E^*$, we have $\sum |\phi(x_i)| < \infty$; equivalently, if
$$
\sup_n\sup_{|\epsilon_i|\leq
1} \left\| \sum_{i=1}^n \epsilon_i x_i \right\|
< \infty \, .
$$
A series is u.c.\ if any subseries is norm convergent.

\section{Unconditionally converging polynomials}
\label{PUC}

In this Section, we prove that whenever $\LEF = \UC (E,F)$, we also have
the equality $\pkEF = \puckEF$, for all $k\in\N$.

We begin with two lemmas.

\begin{Lemma}
\label{UCc0F}
Assume $\LEF = \UC (E,F)$ and ${\cal L}(E,c_0) = \UC (E,c_0)$.
Then we have ${\cal L}(E,c_0(F)) = \UC (E,c_0(F))$.
\end{Lemma}

\Proof
Suppose $T\in {\cal L}(E,c_0(F))$ is not unconditionally converging, and write
$Tx = (T_nx)_n$, with $T_n\in \LEF$.

We can find a w.u.C. series $\sum x_i$ in $E$ such that $\| Tx_i\| >\delta >0$,
for every $i\in\N$. Then there are $n_i\in \N$ so that $\| T_{n_i}x_i\| >
\delta$ $(i\in\N )$.
Since $T_n\in \UC (E,F)$, we have $\| T_n x_i\| \ra 0$, as $i\ra \infty$,
for all $n\in\N$.
This allows assumption (passing to a subseries) that $(n_i)$ is an
increasing sequence.
Choose $\psi_i\in F^*$, $\|\psi_i\| =1$, with $| \psi_i(T_{n_i}x_i)| >
\delta$ $(i\in \N )$, and define $S \in {\cal L}(c_0(F),c_0)$ by
$S((y_n)): = \left(\psi_i (y_{n_i})\right)$. Then
$$
\| STx_j\| = \sup_i \left| \psi_i(T_{n_i}x_j)\right| \geq \left| \psi_j
(T_{n_j}x_j)\right| > \delta \, .
$$
This implies $ST\not\in \UC (E,c_0)$, a contradiction.
\finesp

      The following result is probably well known.
We include a proof for completeness.

\begin{Lemma}
\label{comcopc0}
A space $E$ contains no complemented copy of $c_0$ if and only if
${\cal L}(E,c_0) = \UC
(E,c_0)$.
\end{Lemma}

\Proof
Suppose $T\in {\cal L}(E,c_0)$ is not unconditionally converging. Then we
can find a subspace $M \subseteq E$ isomorphic to $c_0$ such that $T|_M$, the
restriction of $T$ to $M$, is an isomorphism~\cite[Lemma~1]{Pel3}.
By the separable injectivity of $c_0$, $T(M)$ is complemented in $c_0$.
Letting $S:c_0\ra c_0$ be a projection with $S(c_0) = T(M)$, and defining
$U:E\ra E$ by $U:= \left( T |_M\right) ^{-1}ST$,
we have that $U$ is a projection
with $U(E)=M$.
The converse is clear.
\fin

\begin{Th}
\label{mainPUC}
Whenever $\LEF = \UC (E,F)$, we also have $\pkEF = \puckEF$ for all $k\in\N$.
\end{Th}

\Proof
Suppose first that
$E$ contains a complemented copy of $c_0$. Then $F$ cannot contain
a copy of $c_0$, and so every $F$-valued polynomial is unconditionally
converging \cite[Theorem~2]{GG}.

If $E$ contains no complemented copy of $c_0$, then by Lemma~\ref{comcopc0},
we have ${\cal L}(E,c_0) = \UC (E,c_0)$, and we proceed
by induction on $k$.
Suppose the result is true for $(k-1)$-homogeneous polynomials.
Consider $P\in\pkEF$ and a w.u.C. series $\sum x_n$ in $E$.
By the proof of Lemma~4 in \cite{GG}, it is enough to show that $\| Px_n\|
\ra 0$. We define $T\in {\cal L}\left( E,c_0(F)\right)$ by
$$
Tx:= \left( \hat{P}(x_m,\ldots ,x_m,x)\right) _m \; .
$$
We claim that $T$ is well defined. Indeed, for $x\in E$ fixed, we can give a
polynomial $Q\in {\cal P}(^{k-1}\!E,F)$ by
$$
Q(y):= \hat{P}(y,\ldots ,y,x) \hspace{3em} (y\in E) \, .
$$
By the induction hypothesis, $Q\in {\cal P}_{uc}(^{k-1}\!E,F)$.
In particular, $\| Qx_m\| \ra 0$, and the claim is proved.

By Lemma~\ref{UCc0F}, $T\in \UC \left( E,c_0(F)\right)$. Hence,
$$
\| Px_n\| = \left\|\hat{P}(x_n,\ldots ,x_n)\right\| \leq
\sup_m \left\|\hat{P}(x_m,\ldots ,x_m,x_n)\right\| = \|Tx_n\| \lra 0 \, ,
$$
and the proof is complete.
\finesp

      Recall that $E$ has the {\it hereditary\/} DPP if any closed
subspace of $E$ has the DPP.
A polynomial is {\it completely
continuous at
the origin\/} if it takes weakly null sequences into norm null sequences.

\begin{Cor}
Assume $E$ has the hereditary {\rm DPP}, and $\LEF = \UC (E,F)$.
Then $\pkEF = \pcckEF$ for all $k\in\N$.
\end{Cor}

\Proof
By Theorem \ref{mainPUC}, we have $\pkEF = \puckEF$.
Since $E$ has the hereditary DPP, every unconditionally converging
polynomial on $E$ is
completely continuous at the origin \cite[Proposition~20]{GG}.
Hence, every polynomial on $E$ is completely continuous at the origin.

Let now $(x_n)\subset E$ be a sequence weakly converging to $x$, and
$P\in\pkEF$. Then
$$
P(x_n) = \sum_{i=1}^k \hat{P}(x_n-x)^i(x)^{k-i} + P(x) \, .
$$
Since all the polynomials on $E$ are completely continuous at 0, we conclude
that $P(x_n)\ra P(x)$.
\finesp

This Corollary contains a result of \cite{GJ} stating that whenever
$E$ has the
hereditary DPP, and $F$ contains no copy of $c_0$, we have $\pkEF = \pcckEF$.

\section{Completely continuous polynomials}
\label{PCC}

      The fact that $\LEF = \WCo (E,F)$ does not imply that every polynomial
from $E$ into $F$ be weakly compact.
A simple example is the polynomial $Q\in {\cal P}(^k\! \ell_2,\ell_1)$ given in
Section~\ref{Intr}.
This example also shows that if we have $\LEF = \CC (E,F)$, we need not
have $\pkEF = \pcckEF$ either. Ryan \cite{Ry} proved that if $E$ has
the DPP, then every weakly compact polynomial on $E$ is completely continuous.
A modification of his argument allows us to prove that if $E$ has the DPP,
and $\LEF = \WCo (E,F)$, then we have $\pkEF = \pcckEF$ for all $k$.

We need a previous lemma.

\begin{Lemma}
{\rm \cite[Lemma~1.2]{Ry}}
\label{lemaRyan}
An operator $T:E\ra c_0(F)$, with $Tx=\left( T_n(x)\right) _n$,
is weakly compact if and only if the following two conditions are
satisfied:

{\rm (a)} for every $n$, the operator $T_n:E\ra F$ is weakly compact;

{\rm (b)} for every $x^{**}\in E^{**}$, $\lim_n \| T_n^{**}(x^{**})\| =0$.
\end{Lemma}

      We can now state the main result of the Section.

\begin{Th}
\label{main}
Suppose $E$ has the {\rm DPP}, and $\LEF = \WCo (E,F)$.
Given $k\in\N$ and $A\in \LkEF$, let $(x_1^n),\ldots ,(x_k^n) \subset E$ be weak
Cauchy sequences.
Then the sequence $\left( A(x_1^n, \ldots ,x_k^n)\right) _n$ is norm convergent.
\end{Th}

\Proof
By induction on $k$.
For $k=1$, the result is clear.
Assume it is true for $(k-1)$-linear mappings, and take
$A\in \LkEF$, and weak Cauchy sequences
$\left( x_1^n\right) ,\ldots , \left( x_k^n\right) \subset E$.

We suppose first that one of the sequences is weakly null.
To fix notation, let $x_1^n\ra 0$ weakly. For every $z\in E$, the mapping
\be
E\times \stackrel{(k-1)}{\cdots}\times E & \lra & F \\
(x_1, \ldots ,x_{k-1}) & \longmapsto & A(x_1,\ldots ,x_{k-1},z)
\ee
is $(k-1)$-linear.
By the induction hypothesis,
the sequence $\left( A\left(x_1^n,\ldots ,x_{k-1}^n,z\right) \right) _n$ is
norm convergent.
By \cite[Theorem~2.3 and Lemma~2.4]{AHV}, its limit is $0$.
For completeness, we give a short proof of this fact, valid in our case:
Since $E$ has the DPP, for each $\psi\in F^*$, the
$(k-2)$-linear mapping from $E^{k-2}$ into $E^*$ given by
$$
(x_2,\ldots ,x_{k-1}) \longmapsto \psi\circ A( \,\cdot\, ,x_2,\ldots ,x_{k-1},z)
$$
takes the weak Cauchy sequences $\left( x_2^n\right) , \ldots ,
\left( x_{k-1}^n\right) \subset E$ into the weak Cauchy sequence
$\left( \psi\circ A \left( \,\cdot\, ,x_2^n, \ldots ,x_{k-1}^n,z\right)\right)
_n
\subset E^*$.
Again by the DPP of $E$, since $(x_1^n)$ is weakly null, we have
$\psi\circ A\left( x_1^n,\ldots , x_{k-1}^n,z\right) \ra 0$, i.e.\
the sequence $\left( A\left( x_1^n, \ldots ,x_{k-1}^n,z\right)\right) _n$ is
weakly null.
Since it is norm convergent, the limit must be $0$.

Now, we can define the operator $T:E\ra c_0(F)$ by
$$
Tz:= \left( A\left(
      x_1^n,\ldots , x_{k-1}^n,z\right)\right) _n \mbox{ for each }
      z\in E.
$$
We claim that $T$ is weakly compact. Since each coordinate operator
$T_n:E\ra F$ is weakly compact, it is enough (Lemma~\ref{lemaRyan}) to
show that, for every
$z^{**}\in E^{**}$, we have
$\lim_n \left\| T_n^{**}\left( z^{**}\right)\right\| = 0$.

Consider the mapping
$$
\overline{A}: E\times \stackrel{(k-1)}{\cdots}\times E \lra \LEF
$$
given by
$$
\overline{A}(x_1,\ldots ,x_{k-1}): = A(x_1,\ldots ,x_{k-1},\, \cdot\, )\, .
$$
Then $T_n =\overline{A}\left( x_1^n,\ldots ,x_{k-1}^n\right)$.
Since $\LEF = \WCo (E,F)$, for each $z^{**}\in
E^{**}$, we can define $S_{z^{**}} \in {\cal L}(^{k-1}\! E,F)$ by
$$
S_{z^{**}}(x_1,\ldots ,x_ {k-1}) : = \left( \overline{A}(x_1, \ldots ,
x_{k-1})\right) ^{**}\left( z^{**}\right) \, .
$$

By the induction hypothesis,
$$
\left\| T_n^{**}\left( z^{**}\right)\right\| =
\left\| S_{z^{**}}\left( x_1^n,\ldots ,x^n_ {k-1}\right) \right\|
\lra 0 \mbox{ , as } n \ra \infty \, ,
$$
and the claim is proved.

By the DPP of $E$, $T$ is completely continuous.
Therefore $\left( Tx_k^n\right) _n$ converges to some $w = (w_i) \in c_0(F)$.
In particular,
$$
\left\| A\left( x_1^n,\ldots ,x_k^n\right) - w_n \right\| \stackrel{n}{\lra} 0\,
.
$$
Since $\| w_n\| \ra 0$, we conclude that
$\left\| A\left( x_1^n,\ldots ,x_k^n\right) \right\| \ra 0$.

For the general case, suppose that $\left( x_1^n\right) ,\ldots ,
\left( x_k^n\right) \subset E$ are weak Cauchy sequences, and choose two
increasing sequences of indices $(r_n)$, $(s_n)$. Then,
\be
\lefteqn{\left\| A\left( x_1^{r_n},\ldots ,x_k^{r_n}\right) -
A\left( x_1^{s_n},\ldots ,x_k^{s_n}\right)\right\| \leq} \\
& & \left\| A\left( x_1^{r_n}-x_1^{s_n}, x_2^{r_n},\ldots , x_k^{r_n}\right)
      \right\| + \left\| A\left( x_1^{s_n}, x_2^{r_n}-x_2^{s_n},\ldots ,
      x_k^{r_n}\right) \right\| \\
& & + \cdots + \left\| A\left( x_1^{s_n}, \ldots , x_k^{r_n}-x_k^{s_n}
              \right) \right\|   \\
& & \lra 0 \, .
\ee
Hence, the sequence $\left( A\left( x_1^n,\ldots ,x_k^n\right)\right) _n$
is norm convergent.
\fin

\begin{Cor}
\label{Corpol}
Suppose $E$ has the {\rm DPP} and $\LEF = \WCo (E,F)$.
Then we have $\pkEF = \pcckEF$ for all $k\in\N$.
\end{Cor}

The Theorem and Corollary hold, for instance, in the following cases:

(a) $E=C(K)$ with $K$ stonean (e.g.\ $E = \ell_\infty$), and $F \not\supset
\ell_\infty$.

(b) $E = C(K)$ and $F\not\supset c_0$.

(c) $E^*$ has the Schur property, and $F^* \not\supset \ell_1$.

(d) $E^*$ has the Schur property, and $F$ is weakly sequentially complete.

\section{Extension to the bidual}
\label{exten}

We prove that whenever $E^*$ has the DPP property, and $\LEF = \WCo (E,F)$,
then every polynomial from $E$ to $F$ has an extension to a completely
continuous polynomial from $E^{**}$ to $F$.

The following result will be needed:

\begin{Prop}
\label{dualDPP}
{\rm \cite{CG}}
The dual space $E^*$ has the {\rm DPP} if and only if for every $F$ and $T\in
\WCo
(E,F)$, the second adjoint $T^{**}$ is completely continuous.
\end{Prop}

In the next Theorem, we use the same symbol for a multilinear mapping and
its extensions.

\begin{Th}
\label{ext}
Suppose $E^*$ has the {\rm DPP}, and $\LEF = \WCo (E,F)$.
Then each polynomial $P\in\pkEF$ has an extension $\tilde{P}\in \pcck
E^{**},F)$, with $\|\tilde{P}\| = \| P\|$.
\end{Th}

\Proof
Let $A$ be the symmetric $k$-linear mapping associated to $P$.
We extend $A$ to $E^{**}$ coordinatewise by the Davie-Gamelin procedure
\cite{DG}:
for each fixed $j$, $1\leq j \leq k$, and for each fixed $x_1,\ldots ,
x_{j-1} \in E$, and $z_{j+1}, \ldots ,z_k \in E^{**}$, the operator
$$
x\in E \longmapsto A(x_1,\ldots ,x_{j-1},x,z_{j+1},\ldots ,z_k)
$$
is extended to $E^{**}$ by taking its second adjoint. Since
$\LEF = \WCo (E,F)$, it is clear that the extensions have range in $F$.

We define $\tilde{P}(z): = A(z,\ldots ,z)$, for $z\in E^{**}$.
Easily, for each $\psi\in F^*$, $\psi\circ \tilde{P} = \widetilde{\psi\circ P}$.
Since $\|\widetilde{\psi\circ P}\| = \|\psi\circ P\|$ \cite{DG}, we obtain
$\|\tilde{P}\| = \| P\|$.

Proceeding by induction on $k$, we prove that the extension $A$ is
completely continuous. For $k=1$, the result holds by
Proposition~\ref{dualDPP}.
Assume it is true for the $(k-1)$-linear mappings.
Take weak Cauchy sequences $(z_1^n),\ldots ,(z_k^n) \subset E^{**}$, and
$P\in \pkEF$ with associated $A$ as above. Suppose that
one of the sequences, say $(z_k^n)$ to fix notation, is weakly null.
For every $z\in E^{**}$, the mapping
$$
(z_2,\ldots ,z_k) \in (E^{**})^{k-1} \longmapsto A(z,z_2,\ldots ,z_k)
$$
is completely continuous, by the induction hypothesis.
Therefore, the sequence
$\left( A(z,z_2^n,\ldots ,z_k^n)\right) _n$ converges to zero, as in
Theorem~\ref{main}.

Defining $T:E\ra c_0(F)$ by $Tx:= \left(A(x,z_2^n, \ldots , z_k^n)\right)_n$,
for $x\in E$,
we have
$$
\| T_n^{**}(z)\| = \| A(z,z_2^n,\ldots ,z_k^n)\| \lra 0 \, .
$$
By Lemma~\ref{lemaRyan}, $T$ is weakly compact.
By Proposition \ref{dualDPP}, $T^{**}$ is completely continuous.
Hence, as in Theorem~\ref{main}, $A(z_1^n,\ldots ,z_k^n) \ra 0$.
The proof finishes as in Theorem~\ref{main}.
\fin

\end{document}